\documentclass[11pt]{article}
\usepackage{amsmath,amsthm,amsfonts,amssymb,amscd,amsbsy}
\usepackage[utf8]{inputenc}
\usepackage[pdftex]{graphicx}
\usepackage{setspace}
\usepackage[bottom]{footmisc}
\usepackage{gensymb}
\usepackage{lipsum}
\usepackage{latexsym}
\usepackage{fancyhdr}
\everymath=\expandafter{\the\everymath\displaystyle}

\newcommand{\R}{\mathbb{R}}

\newcommand{\Fi}{\Phi}
\newcommand{\abre}{\left\langle}
\newcommand{\fecha}{\right\rangle}

\newcommand\blfootnote[1]{%
  \begingroup
  \renewcommand\thefootnote{}\footnote{#1}%
  \addtocounter{footnote}{-1}%
  \endgroup
}

\renewcommand*{\thefootnote}{\fnsymbol{footnote}}

\title{Convex geometric reasoning for crystalline energies}

\author{Thaicia Stona\\\\ University at Buffalo - State University of New York\\\\ thaicias@buffalo.edu}

\date{ \blfootnote{This work was initiated and mostly developed during the Thematic Program on Variational Problems at the Fields Institute, Toronto, Canada (Fall 2014). The author is deeply grateful to the Fields Institute for the support received and its hospitality. The author would like to thank Almut Burchard, Jean Ellen Taylor, Robert McCann and all the supporters of the Maths, Metallurgy \& Crystals project (please refer to the complete list at http://goo.gl/GTgI1Y).}  }

\begin{document}
\maketitle
\begin{abstract}

$ $

\noindent ``The present work revisits the classical Wulff problem restricted to crystalline integrands, a class of surface energies that gives rise to finitely faceted crystals. The general proof of the Wulff theorem was given by J.E. Taylor (1978) by methods of Geometric Measure Theory. This work follows a simpler and direct way through Minkowski Theory by taking advantage of the convex properties of the considered Wulff shapes.'' \ \textit{(Final version is published on Caspian Journal of Computational \& Mathematical Engineering, 1, 2016)}

$ $

\end{abstract}

\begin{center}
\textbf{Introduction}
\end{center}

\noindent This work is a short though sufficiently self-contained incursion into the Wulff construction and the Wulff theorem for faceted crystals, mathematically represented by the class of crystalline integrands.\\

\noindent The aim of the Wulff Problem is to find a surface whose total surface energy is minimal for a given fixed volume. This classical problem is also known as the Equilibrium Shape Problem, and the solution is also called an equilibrium shape, or simply a crystal. The problem is named after George Wulff, who invented an algorithm to determine the final shape of a crystal that grows near equilibrium, based on Josiah Willard Gibbs principle of the surface Gibbs free energy minimization for the evolution of a crystal droplet.\\

\noindent By convexifying $\gamma$, one induces a $\gamma$-metric on the dual space of the solution space. Since such a class of problems has polyhedral solutions, we can dismiss the geometric measure versions of Brunn-Minkowski Theorem from Federer and Wulff Theorem by applying the Legendre transform to the canonical version of the Wulff construction and build our way to the Convex Geometry version of Brunn-Minkowski Theorem through geometric inequalities and convexity. We show equivalences between constructions and some relations between the crystalline integrand and the area integrand version of the problem - isoperimetry and minimal surfaces.\\

\begin{center}
\textbf{I. The Crystalline Variational Problem}
\end{center}

$ $

\noindent The Wulff shape arises in surface energy minimization problems when the energy function $\Phi$ is anisotropic. For isotropic energies and a given amount of mass, the equilibrium shape is well known: a ball, formally denoted by the n-dimensional sphere $S^{n-1}$. This shape encloses the prescribed mass whereas minimizing the surface area of it. This is stated by the classic Isoperimetric Inequality.\\

\noindent For anisotropic energies, the analogous minimizer is the Wulff shape. Such energies have been heuristically misrepresented by simple functions that very often are not well-defined, presenting many singularities and unbounded energy spots. We avoid this imposture following Taylor's Geometric Measure Theory characterization of the energy. In this context, the surface energy, also called the energy function of the anisotropic problem, is an integrand, as defined below:\\

\noindent \textbf{Def.:}[integrand] An integrand on $\R^{n+1}$ is a function that will represent the surface energy function
\begin{center}
$\Fi : \R^{n+1}\times G_0(n+1,n) \longrightarrow [0,+\infty]$\\
\end{center}

$ $

\noindent where the Grassmannian $G_0(n+1,n)$ is the manifold that parametrizes every n-dimensional linear subspace of $\R^{n+1}$, ie, all the hyperplanes of the (n+1)-dimensional Euclidean space.\\

\noindent An integrand is defined constant coefficient iff $\Fi(x,\pi) =  \Fi(p,\pi)$ $\forall x, p \in \R^{n+1}$, $\forall \pi \in G_0(n+1,n)$. In this case, $\Fi$ is a function of its second variable only. An integrand is unoriented if it independs on the orientation of $\pi$. We will assume all integrands are continuous, constant coefficient and positively oriented.\\

\noindent \textbf{Def.:}[Wulff construction] Given an integrand $\Fi$, plot it radially by taking each direction $v\in S^n$ and calculating $\Fi$ on the positively oriented plane $\pi$ whose normal vector is $v$, $\pi = \left\{ x\in \R^{n+1} / \left\langle x , v \right\rangle = 0 \right\}$. We will denote $\pi$ by $v^\bot$ and vice versa. Plot $\Fi(v^\bot)$ in $v$ direction: $\Fi(v^\bot) \ v$. Then, for each $v$, define the half-space $H_v \doteq \left\{ x\in \R^{n+1} / \left\langle x , v \right\rangle \leq \Fi(v^\bot) \right\}$.\\

\noindent Take the intersection of all half-spaces. The resulting set $W_\Fi$ is the Wulff shape of $\Fi$, also called the crystal of $\Fi$:

$$
W_\Fi \doteq \bigcap_{v \in S^n} H_v
$$

\noindent For an isotropic energy, $\Fi \equiv$ constant: the Wulff problem reduces to the Isoperimetric Inequality and the crystal is an Euclidean ball; that is the case of a soap bubble.\\

\noindent \textbf{Obs.:} We can extend homogeneously the function $\Fi$ in order to calculate it on other planes related to non unitary direction vectors by formalizing the explained abuse of notation defining the dual function $\Fi^\star$ as follows

$$
\Fi^\star  : S^n \longrightarrow [0,+\infty] ,  \ \ \ \ \Fi^\star(v) = \Fi(\pi)
$$

$ $

\noindent where $v \bot \pi$ as defined above, $\Fi^\star (p) \doteq |p| \ \Fi^\star(\frac{p}{|p|})$.\\

\noindent Note that since $W_\Fi$ is given by an intersection of half-spaces, then $W_\Fi$ is convex. Also we can assume $0\in W_\Fi$ always. The physical meaning of the origin is the crystal seed for growing a crystal, a tiny monocrystal that induces the orientation of the new crystal.\\

\noindent \textbf{Def.:}[Legendre Transform] Let $\xi : S^{n-1} \longrightarrow \R^+$ be a continuous function. The (first) Legendre transform of $\xi$ is

$$
\xi^\star(v) \doteq \inf_{\abre \theta, v \fecha >0} \frac{\xi(\theta)}{\abre \theta, v \fecha} \ \ \mbox{ where } \ |\theta| = 1
$$

$ $

\noindent An alternative construction of the Wulff shape is based on the Legendre transform as in [5]:\\

\noindent \textbf{Def.:}[Fu's Wulff construction] Let $W$ be the operator over integrands
$$
W(\Fi)(\pi) \doteq \inf_{v\in S^n} \frac{\Fi^\star(v)}{\abre \pi^\bot , v\fecha}
$$

\noindent where $\abre \pi^\bot , v\fecha > 0$. Then the crystal of $\Fi$ is the set enclosed by the radial plot of $W(\Fi)$, plotted as explained before. Also, the orientation of $W_\Fi$ is defined positive.\\

\noindent \textbf{Proposition:} The two given definitions of crystal are equivalent.\\

\noindent \textbf{Proof:} Call $Z$ the operator defined in Fu's construction instead of $W$:\\

\noindent $(W_\Fi \subset Z_\Fi ):$ Let $y\in W_\Fi$, ie, $\abre y, v \fecha \leq \Fi^\star(v)$  \ ($\forall v\in S^n$). Then

$$
\abre y, v \fecha =\frac{|y|}{|y|} \abre y, v \fecha = |y| \abre \frac{y}{|y|}, v \fecha \leq \Fi^\star(v)
$$

$ $

\noindent If $\abre \frac{y}{|y|}, v \fecha > 0$ then we have $  |y| \leq \frac{\Fi^\star(v)}{\abre \frac{y}{|y|}, v \fecha}$. Since the inequality holds for arbitrary $v$, then

$$
|y| \leq \inf_{v\in S^n} \frac{\Fi^\star(v)}{\abre y , v\fecha} \mbox{ \ , ie, \ }y\in Z_\Fi
$$

\noindent If $\abre \frac{y}{|y|}, v \fecha \ngeq 0$, then obviously the inequality holds, with $y$ in the same half-space bounded by $\abre x , y \fecha = \Fi^\star (\frac{y}{|y|})$;\\

\noindent $(Z_\Fi \subset W_\Fi ):$ Let $y\in Z_\Fi$, ie, $|y|\leq (Z(\Fi))^\star (\frac{y}{|y|})$, $\abre \frac{y}{|y|}, v \fecha > 0$. Then:
$$
|y| \leq \frac{\Fi^\star(v)}{\abre \frac{y}{|y|}, v \fecha} \ \ \ \forall v \in S^n \iff |y|\abre \frac{y}{|y|}, v \fecha = \abre y, v \fecha \leq \Fi^\star (v)
$$

\noindent for all $v\in S^n$; but then $y\in H_v \ \ (\forall v\in S^n) \ \ \Rightarrow \ \  y\in W_\Fi \ \ \ \Box$

$ $


\begin{center}
\textbf{II. Pathway through convexity}
\end{center}

$ $

\noindent It is easy to visualize what kind of Wulff shape one gets when the intersection of half-spaces is finite: a polyhedron, except for unbounded and/or empty intersections. That is the case of anisotropic energies: we say that an integrand $\Fi$ is crystalline if its Wulff shape, or crystal $W_\Fi$ is a polyhedron. Now we take advantage of this fact:\\

\noindent \textbf{Def.:}[extreme point] Given a set $K\subset \R^n$, $x\in K$ is extreme if it cannot be expressed as a convex combination of any two other points of $K$.\\

\noindent \textbf{Def.:}[polytope] A polytope $P\subset \R^n$ is the convex hull of a finite set:  $P = [\left\{p_1,p_2,...,p_k \right\}]$.\\

\noindent \textbf{Def.:}[polar body] Given a convex set $K$, the polar body of $K$ is the set $K^\star \doteq \left\{ x\in\R^n / \abre x,y \fecha \leq 1 \ (\forall y \in K ) \right\}$.\\

\noindent \textbf{Lemma:}  A supporting hyperplane $H$ to a bounded convex set $K$ contains at least one extreme point of $K$.\\

\noindent \textbf{Proof:} Denote the set of extreme points of $K$ by $E_K$. Since $K$ is convex, $K=[K]$, so   $E_k\subset K \Rightarrow [E_K]\subset K $. We also have that $H\cap K = H\cap \partial K$, so the set of extreme points of $H\cap K$, $E_{H\cap K}$, is the set $E_{H\cap E_K}$. Now suppose the claim is true for every set with dimension $\leq m-1$. Then it is also true for all sets of dimension $m$, since if a given non-extreme point in $m$ dimension could be written as a convex combination in dimension $m-1$, then it would be sufficient to write it in $m$ dimension putting $\lambda_m = 0$. But for dimension 1, the claim is trivially true. Therefore it is true for any dimension. \ \ \ $\Box$

\noindent \textbf{Theorem 1:}  A bounded convex set $K$ is the convex hull of its extreme points.\\

\noindent \textbf{Proof:} Since $E_k\subset K \Rightarrow [E_K]\subset K $, we only need to prove that $K\subset [E_K]$. Suppose some $x\in K$ is not in $[E_K]$. Then there exists a separating hyperplane H that separates strictly $x$ from $E_K$. The parallel supporting hyperplane of $K$ that is strictly separated from $E_K$ by $H$ must contain a point of $E_K$ (lemma). Contradiction.  \ \ \ $\Box$

\noindent \textbf{Corollary:} Every polytope is a finite intersection of half-spaces.\\

\noindent \textbf{Proof:} If $P$ is finite, then so is $E_P$. For each $p\in E_P$, let $A_p$ be the set of supporting hyperplanes that contains $p$ and also contains at least another extreme point of $P$. Then take $A_p'\subset A_p$ the subset that contains supp. hyperplanes intersecting the maximum number of extreme points as possible (this number is well-defined since the very $\#P$ is a majorant). The facets of $P$ will be contained on those hyperplanes; for each facet define the half-space oriented to contain the origin and take the intersection of it. Because of the theorem, $P$ is contained in this intersection.  \ \ \ $\Box$

\noindent \textbf{Theorem 2:} If $K$ is convex, then $K^{\star\star} = K$\\

\noindent \textbf{Proof:}\\
\noindent $(K \subset K^{\star\star} )$ \ Let $x\in K$. Then for any $y\in K^\star$ we have $\abre x,y \fecha \leq 1$. But then, since $x$ is arbitrary, it has to be in $K^{\star\star}$.\\

\noindent $(K^{\star\star} \subset K )$ \ Let $y\in K^{\star\star}$ and suppose $y\notin K$. Then there is a separating hyperplane $H$ that separates $y$ from $K$, $H=\left\{ x \ / \abre x,v \fecha = 1  \right\}$,
$$ \abre x,v \fecha \leq 1 \mbox{ when } x\in K \ \ \mbox{ and } \ \ \abre y,v \fecha > 1$$

\noindent But if $ \abre x,v \fecha \leq 1$ when  $x\in K$, then $v\in K^\star$ and $\abre y,v \fecha \leq 1$ since $y\in K^{\star\star}$. Contradiction.  \ \ \ $\Box$

\noindent Theorem 2 reveals a link between Convex Geometry and Functional Analysis: given a polyhedral crystal $W$, we apply the corollary to define a convex $\Fi_C$ whose crystal coincides with $W$, so that $\Fi_C$ is the "smallest" enclosing function for $W$. For that, we use the theorem 2 by taking the polar of $W$. Since $\Fi$ is a linear operator, we know its behavior everywhere by homogeneous extension. By Riesz representation theorem, the crystal $W$ is the polar of the unit ball $\Fi_C \equiv 1$ in the dual norm. That gives us the surface energy scaled so that the Wulff shape is given in units of surface free energy.\\

\noindent \textbf{Def.:}[Steiner symmetrization] For a convex body $K \subset \R^n$ and a $\theta\in S^{n-1}$, the Steiner symmetrization of $K$ in the direction of $\theta$ is given by

$$
S_\theta (K) \doteq \left\{  x+\lambda.\theta \ | \ x\in Proj_{\theta^\perp} K, \ \lambda \in \R  \right\}
$$

\noindent where $|\lambda| \leq \frac12 |K \cap \left\{  x+ \R\theta  \right\}|$. Some properties are the fact that $|S_\theta (K)|=|K|$, $S_\theta(K)$ is convex and the convex Minkowski sum of symmetrizations equals to the symmetrization of the convex sum of the bodies. The symmetrization process slices $K$ along $\theta$, aligning the slices by putting their midpoints in $\theta^\perp$.\\

\setcounter{footnote}{1}

\begin{center}
\includegraphics[scale=0.6]{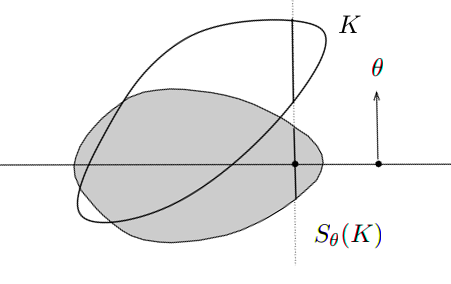}

$ $

\noindent Figure 1: Example of Steiner symmetrization of $K$ along the vector $\theta$ \footnote{Figure adapted from [24]}\\

\end{center}

\noindent A useful classical result is stated below without its proof, which follows directly from the several interesting properties of the Steiner Symmetrization process. A more curious reader might refer to [6], [20] or [21].\\

\noindent \textbf{Theorem:}[Steiner-Schwarz] Given a convex body $K\subset\R^n$ and $F$ a\\ k-dimensional subspace, then there exists a sequence of symmetrizations $\theta_j$ such that the limiting body $\overline{K}$ satisfies $|\overline{K} \cap \left\{  x+ F  \right\}| = |K \cap \left\{  x+ F  \right\}|$, where $\overline{K} \cap \left\{  x+ F  \right\}$ is a k-dimensional ball centered in $x$ with radius $r(x)$.\\

\noindent \textbf{Theorem:}[Brunn's Concavity Principle] Given $K \subset \R^n$ a convex body  and $F$ a k-dimensional subspace of $\R^n$, the function $f: F^\perp \longrightarrow \R^+ $ given by $f(x) = |K \cap \left\{  x+ F  \right\}|^\frac1n$ is concave on its support.\\

\noindent \textbf{Proof:} Apply the former theorem and use that $supt \ r(x) = Proj_{F^\perp} K$, \ $f(x) = |\overline{K} \cap \left\{  x+ F  \right\}| = Vol(S^{k-1}) = \frac{\pi^\frac{k}{2}}{\Gamma(\frac{k}{2} +1)}r(x)^k$.  \ \ \ $\Box$

\begin{center}
\includegraphics[scale=0.5]{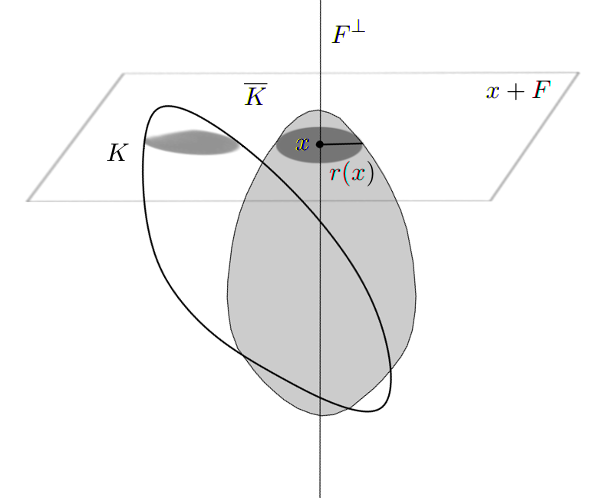}

$ $

\noindent Figure 2: Application of Steiner-Schwarz to prove Brunn's Concavity Principle, where $n=3$, $k=2$ \footnote{Figure adapted from [24]}\\

\end{center}

\noindent The Brunn-Minkowski inequality is the crucial ingredient for proving the optimality of the Wulff shape. We conclude this section with a proof based on convex sum of two convex bodies and the Concavity Principle:\\

\noindent \textbf{Theorem:}[Brunn-Minkowski Inequality] Given non-empty compact subsets $A, B$ of $\R^n$

$$
|A+B|^\frac1n \geq |A|^\frac1n + |B|^\frac1n
$$

$ $

\noindent \textbf{Proof:} Take the Steiner symmetrization of $A$ and $B$ to find two convex bodies in $\R^n$. Create their convex sum L on $\R^{n+1}$ by taking the convex hull of $S_\theta(A) \times {0}$ and $S_\theta(B) \times {1}$, where $0, 1$ belong to the additional real axis for the convex sum, so that $L(t) = \left\{ x\in\R^n | (x,t) \in L \right\}$.\\ Then $L(\frac12) = \frac{S_\theta(A)}{2} + \frac{S_\theta(B)}{2} = \frac{S_\theta(A) + S_\theta(B)}{2}$. By the concavity principle applied for $F=\R^n$

$$
\left| \frac{S_\theta(A) + S_\theta(B)}{2} \right|^\frac1n \geq \frac12\left|S_\theta(A)\right|^\frac1n + \frac12\left|S_\theta(B)\right|^\frac1n
$$
$ $  \ \ \ $\Box$\\

\begin{center}
\textbf{III. The Wulff Theorem}
\end{center}

$ $

\noindent Wulff's 1901 seminal article provided a method to predict crystal shapes after Gibbs' proposition on the minimization of surface energy; since then, many have worked on the subject. Nevertheless, it was Taylor ([1]) who proved that the Wulff construction determines the unique minimizer $W_\Fi$ for the integral of $\Fi$ over the boundary $\partial W_\Fi$. The proof requires some concepts from Geometric Measure Theory, which are now introduced:\\

\noindent \textbf{Def.:}[integral current] An integral n-dimensional current $S\subset \R^{n+1}$ is a rectifiable oriented hypersurface generalized through GMT so that eventual anomalous portions are still well-behaved enough to allow integration with respect to the measure $|S|$ on $\R^{n+1}$, which is a function of the Hausdorff measure $H^n$ restricted to the support of $S$, which can be arbitrarily closely approximate by a n-d $C^1$ manifold. An interesting property of currents is that their boundaries also have the essential properties to allow boundary integration (for more see [3]). In the next theorem $T$ will denote the current whose boundary is an integral current. The total surface energy of an integral current $S\subset \R^{n+1}$ is given by:

$$
\Fi (S) \doteq \int_{x\in S} \Fi[n_S(x) ] dH^nx
$$

$ $

\noindent We also define for $h>0$ the homothety in $\R^{n+1}$ $\mu_h(x) = hx$ and the integrand $\Fi$ the isomorphism $W^h_\Fi = \mu_{h\sharp} (W_\Fi)$ following [1].\\

\noindent \textbf{Theorem:}[Wulff] Given an integrand $\Fi$, then for every n-dimensional current $P \subset \R^{n+1}$
$$
\Fi(\partial W_\Fi) \leq \Fi(\partial P)
$$

\noindent up to translations and homotheties, such that their mass coincide, $M(P)=M(W_\Fi)$\\

\noindent \textbf{Proof:} Let $P$ be a current with $\partial P$ its positively oriented, piecewise $C^1$ boundary. Then

$$\Fi(\partial P) = \int \Fi (\overrightarrow{\partial P}(x)) d|\partial P|x \geq \int supt \ (W_\Fi) (\overrightarrow{\partial P}(x)) d|\partial P|x$$
$$ = \lim_{h\rightarrow 0} \frac{M(P^h)-M(P)}{h}$$

\noindent where  $M(W_\Fi)=M(P) \mbox{ and } M(W^h_\Fi)=h^{n+1}M(W_\Fi)$ and $P^h$ is the positively oriented current given by the Minkowski sum $x+y$ where $x\in supt  \ P$ and $ y\in supt \ W^h_\Fi$. Then Brunn-Minkowski inequality implies:

$$ = \lim_{h\rightarrow 0} \frac{M(P^h)-M(P)}{h} \geq  \lim_{h\rightarrow 0} \frac{(1+h)^{n+1}M(W_\Fi)-M(W_\Fi)}{h}
$$
$$ = \lim_{h\rightarrow 0} \frac{(1+h^{n+1}-1)}{h}M(W_\Fi) = \lim_{h\rightarrow 0} \frac{M(W_\Fi)}{h}\sum_{i=1}^{n+1} h^i.\frac{(n+1)!}{i!(n+1-i)!}$$
$$=\lim_{h\rightarrow 0} M(W_\Fi).(n+1) = (n+1).M(W_\Fi) $$

$ $

\noindent In particular for $P=W_\Fi$, the above inequalities are equalities. By using the fact that $M(P) = M(W_\Fi)$ , we conclude that $ \Fi (\partial W_\Fi) \leq \Fi (\partial P)$.

\noindent Such shape is unique modulo translations and homotheties, and since the mass is fixed, follows the uniqueness of $W_\Fi$.  \ \ \ $\Box$

$ $

\begin{center}
\textbf{IV. Conclusion}
\end{center}

$ $

\indent In this exposition, different fundamental areas of Mathematics were gathered to structure a simple mathematical basis for the equilibrium shape problem with a crystalline integrand. A natural generalization of the Wulff construction for non-equilibrium growth is to replace the energy function for the correspondent potential that controls the process, the mobility function. Also, through Kinectic PDEs, a flourishing area of mathematical modelling in the Sciences, it might be of interest to study the growth and the stability of such shapes.

\begin{center}
\textbf{V. References}
\end{center}

\begin{enumerate}
\item Taylor, J.E. \textit{Crystalline variational problems}, 1978
\item Burchard, A. \textit{A short course on rearrangement inequalities}, 2009
\item Federer, H. \textit{Geometric Measure Theory}, 1969
\item McCann, R. \textit{Equilibrium shapes for plannar crystals in an external field}, 1998
\item Fu, J. \textit{A mathematical model for crystal growth and related problems}, 1976
\item Brazitikos, S., Giannopoulos, A., Valettas, P., Vritsiou, B. \textit{Geometry of Isotropic Convex Bodies}, 2014
\item Gibbs, J.W. \textit{Collected Works Vol.1}, 1948
\item Wulff, G. \textit{Zeitschrift fur Krystallographie und Mineralogie}, 1901
\item Taylor, J.E., Cahn, J.W., Handwerker, C.A. \textit{Evolving crystal forms: Frank's characteristics revisited}, 1991
\item Wills, J.M. \textit{Wulff-Shape, Minimal Energy and Maximal Density}, 2001
\item Micheletti, A., Patti, S., Villa, E. \textit{Crystal Growth Simulations: a new Mathematical Model based on the Minkowski Sum of Sets}, 2005
\item Taylor, J.E. \textit{Crystalline Variational Methods}, 2002
\item Cahn, J.W., Handwerker, C.A. \textit{Equilibrium geometries of anisotropic surfaces and interfaces}, 1993
\item Cahn, J.W., Hoffman, D.W. \textit{A vector thermodynamics for anisotropic surfaces - II. curved and faceted surfaces}, 1974
\item Palmer, B. \textit{Stable closed equilibria for anisotropic surface energies: Surfaces with edges}, 2011
\item Koiso, M., Palmer, B. \textit{Stable surfaces with constant anisotropic mean curvature and circular boundary}, 2013
\item Craig Carter, W., Taylor, J.E., Cahn, J.W. \textit{Variational Methods for Microstructural Evolution}, 1997
\item Herring, C. \textit{Some theorems on the free energies of crystal surfaces}, 1951
\item Almgren, F., Taylor, J.E., Wang, L. \textit{Curvature driven flows: a variational approach}, 1993
\item Eggleston, H.G. \textit{Convexity}, 1958
\item Schneider, R. \textit{Convex Bodies: The Brunn–Minkowski Theory}, 2014
\item Peng, D., Osher, S., Merriman, B., Zhao, H. \textit{The geometry of Wulff Crystals
Shapes and its relations with Riemann problems}, 1998
\item Micheletti, A., Burger, M. \textit{Stochastic and deterministic simulation of nonisothermal crystallization of polymers}, 2001
\item Burchard, A. \textit{How to achieve radial symmetry through simple rearrangements}, 2012

\end{enumerate}

\end{document}